\documentclass[twoside, 12pt]{article}
\usepackage{amssymb, amsmath, mathrsfs, amsthm}
\usepackage{graphicx}
\usepackage{color}
\usepackage[top=2cm, bottom=2cm, left=2cm, right=2cm]{geometry}
\usepackage{float}

\newcommand{\bd}{\begin{description}}
\newcommand{\ed}{\end{description}}
\newcommand{\bi}{\begin{itemize}}
\newcommand{\ei}{\end{itemize}}
\newcommand{\be}{\begin{enumerate}}
\newcommand{\ee}{\end{enumerate}}
\newcommand{\beq}{\begin{equation}}
\newcommand{\eeq}{\end{equation}}
\newcommand{\beqs}{\begin{eqnarray*}}
\newcommand{\eeqs}{\end{eqnarray*}}

\definecolor{DarkGreen}{rgb}{0.2, 0.6, 0.3}

\newtheorem{theorem}{Theorem}
\newtheorem{conjecture}{Conjecture}

\newtheorem{lemma}{Lemma}

\newtheorem{case}{Case}

\begin{document}
\begin{center}
    {\Large \bf Connectivity keeping paths for $k$-connected bipartite graphs}\\

    \vspace{10mm}

    {\large \bf Meng Ji}

    \vspace{5mm}

    \baselineskip=0.20in

    {\it College of Mathematical Science, Tianjin Normal University, \\
       Tianjin, China}\\
    {\rm E-mail:} {\tt mji@tjnu.edu.cn }\\[2mm]
    
\vspace{5mm} 
 \end{center}

\vspace{5mm}

\baselineskip=0.23in

\begin{abstract}
Luo, Tian and Wu [Discrete Math. 345 (4) (2022) 112788]  conjectured that for any tree $T$ with bipartition $(X,Y)$,
every $k$-connected bipartite graph $G$ with minimum degree at least $k+w$, where $w=\max\{|X|,|Y|\}$,
contains a tree $T'\cong T$ such that $\kappa(G-V(T'))\geq k$. In the paper, we confirm the conjecture when $T$ is an odd path on $m$ vertices. We remind that Yang and Tian \cite{YT2} also prove the same result by a different way.
\\[2mm]
\end{abstract}

\section{Introduction}

In this paper, all the graphs are finite, undirected and simple. For
graph-theoretical terminology and notations not defined here, we
follow \cite{BM}. Some basic symbols and definitions are needed
to be introduced. The minimum degree and the connectivity of
a graph $G$ are denoted by $\delta(G)$ and $\kappa(G)$, respectively. The vertex set and edge set of $G$ are denoted by $V(G)$ and $E(G)$,
respectively. For $V\subseteq V(G)$, we denote the minimum degree of the vertex set $V$ in $G$ by $\delta_{G}(V)$. 
 For a subgraph $H\subseteq G$ and a subset
$V'\subseteq V(G)$, write $V'\cap H$ for $V'\cap V(H)$. For
$H\subseteq G$, set $\delta_{G}(H):=\min_{x\in H}{d_{G}(x)}$.

 Let $S$ be a separating set of $G$, and the cardinality of a minimum separating set $S$ is
denoted by $|S|=\kappa(G)$. If $S$ is a minimum separating set of
the graph $G$, then the union components $F$ of at least one component of $G-S$ with
$G-S-F\neq \emptyset$ is called a fragment to $S$, and the
complementary fragment $G-(S\cup V(F))$ is denoted by $\bar{F}$. For any minimum separating set $S'$ and a fixed minimum separating set $S$,
if a fragment of $G$ to $S$ does not properly contain another fragment of $G$ to $S'$, then
it is called an end of $G$. Clearly, every graph contains an end
except for the complete graphs. For a fragment $F$ of $G$ to $S$, it
definitely follows $S=N_{G}(F)$.

In 1972, Chartrand, Kaigars and Lick \cite{CKL} proved that
\emph{every $k$-connected graph $G$ with $\delta (G)\geq
\lfloor\frac{3k}{2}\rfloor$ has a vertex $x$ with $\kappa(G-x)\geq
k$}. Nearly forty years later, Fujita and Kawarabayashi \cite{SK}
extended the result, that is, every $k$-connected graph $G$ with
$\delta (G)\geq \lfloor\frac{3k}{2}\rfloor+2$ has an edge $xy$ such
that $G-\{x,y\}$ remains $k$-connected. Furthermore, Fujita and
Kawarabayashi \cite{SK} posed a conjecture that
\emph{for all positive integers $k$, $m$, there is a (least) non-negative
integer $f_{k}(m)$ such that every $k$-connected graph $G$ with
$\delta(G)\geq \lfloor\frac{3k}{2}\rfloor+ f_{k}(m)-1$ contains a
connected subgraph $W$ of order $m$ such that $G-V(W)$ is still
$k$-connected.}

In 2010, Mader \cite{mader1} confirmed that conjecture posed by Fujita and
Kawarabayashi and he proposed the following conjecture.

\begin{conjecture}[Mader \cite{mader1}]\label{conj0}
For every positive integer $k$ and finite tree $T$ of order $m$,
every $k$-connected finite graph $G$ with minimum degree
$\delta(G)\geq\lfloor\frac{3k}{2}\rfloor+m-1$ contains a subgraph
$T'\cong T$ such that $G-V(T')$ remains $k$-connected.
\end{conjecture}

For Mader's conjecture, although many scholars have proved plenty of results partially confirming the conjecture, there are few results to partially confirm it for $k\geq 4$.  For more results, the reader can refer to  \cite{AN,H,TK,HL2,Lu,mader1,mader2,Tian1,Tian2}. 
If $G$ is a bipartite graph, Luo, Tian and Wu \cite{LTW} proved that \emph{every $k$-connected bipartite graph $G$ with minimum degree at least $k + m$
contains a path $P$ with order $m$ such that $\kappa(G-V(P))\geq k$.} Based on that, they proposed a general conjecture. 
\begin{conjecture}{\upshape\cite{LTW}}\label{conj}
For any tree $T$ with bipartition $(X,Y)$, every $k$-connected bipartite
graph $G$ with $\delta(G)\geq k+w$, where $w = \max\{|X|,|Y|\}$, contains a tree $T'\cong T$ such that
$$\kappa(G-V(T'))\geq k.$$
\end{conjecture}

The case that $T$ is caterpillar and $k\leq 2$ has been proved by Zhang \cite{Z}. Additionally, Yang and Tian \cite{YT} verified Conjecture \ref{conj} when $T$ is caterpillar and $k = 3$, or $T$ is a spider and $k\leq 3$. Just as the research of Conjecture \ref{conj0}, there are no results on high connectivity graph to confirm Conjecture \ref{conj}. In the paper, the author will verify Conjecture \ref{conj} when $T$ is an odd path. 

\begin{theorem}\label{main-thm}
Let $m$ be an odd positive integer. Every $k$-connected bipartite
graph $G$ with $\delta(G)\geq k+\lceil\frac{m}{2}\rceil$ contains a path $P$ on $m$ vertices such that
$$\kappa(G-V(P))\geq k.$$
\end{theorem}

\section{Proof of Theorem \ref{main-thm}}
We first prove our key lemma as follows.
\begin{lemma}\label{lem1}
For every positive integer $k$,
let $G$ be a bipartite graph with $\kappa(G)=k$ and
$\delta(G)\geq k+1$. Then we have $F\cap S_{1}=\emptyset$ for each end $F$ of $G$ to a minimum separating set $S$ and for any minimum separating set $S_{1}$ in $G$.
\end{lemma}

\begin{proof}
Let $G$ be a bipartite graph with a bipartition $(X',Y')$. 
Let $S_{1}\neq S$ be a minimum vertex cut of $G$ and denote by $F_{1}$ a fragment of
$G$ to $S_{1}$. Clearly, $V(F_{1}),$ $S_{1}$, and $V(\bar{F}_{1})$ (resp., $V(F),$ $S,$ and $V(\bar{F})$) form a partition of $V(G)$. (For convenience, we will abbreviate $V(F_{1})$, $V(\bar{F}_{1})$, $V(F),$ and $V(\bar{F})$ as $F_{1}$, $\bar{F}_{1}$, $F,$ 
and $\bar{F}$, respectively)
The graph $G$ could be partitioned to at most nine parts, namely,
$\{F\cap F_{1}$, $F\cap S_{1}$, $F\cap \bar{F}_{1}$, $S\cap F_{1}$, $S\cap S_{1}$, $S\cap \bar{F}_{1}$,
$\bar{F}\cap F_{1}$, $\bar{F}\cap S_{1}$, $\bar{F}\cap \bar{F}_{1}\}$.
We suppose, to the contrary, that $F\cap S_{1}\neq\emptyset$.
Three cases will be considered.

\begin{case}
 $|S_{1}\cap F|=|\bar{F}\cap S_{1}|$.
\end{case}
To begin with, let $|S\cap F_{1}|=|\bar{F}_{1}\cap S|$ as a constraint. If at most one of $F\cap F_{1}$ and $F\cap \bar{F}_{1}$
is empty-set, then, obviously, $(S\cap F_{1})\cup(S\cap S_{1})\cup(F\cap S_{1})$ is also a minimum separating
set of $G$. Without loss of generality, assume that $F\cap F_{1}\neq \emptyset$. It implies that $F\cap F_{1}$ is a fragment, which contradicts the fact that $F$ is an end.
If $F\cap F_{1}=\emptyset $ and $F\cap \bar{F}_{1}=\emptyset$, then $F\subset S_{1}$.
By $\delta(F)\geq \delta(G)-|S|\geq 1$, $F$ contains at least an edge $xy$. Because of $\delta_{G}(X')\geq k+1$ and $\delta_{G}(Y')\geq k+1$, we have $d_{G}(x)+d_{G}(y)> 2k+1$. Meanwhile, since $x$ and $y$ have distinct neighbors in $S$, we also have $d_{G}(x)+d_{G}(y) \leq |S|+|F|\leq |S|+|S_{1}|\leq 2k$, a contradiction.

Suppose, without loss of generality, that $|S\cap F_{1}|>|\bar{F}_{1}\cap S|$.
If $F\cap \bar{F}_{1}\neq \emptyset$, then $(F\cap S_{1})\cup (S_{1}\cap S)\cup (\bar{F}_{1}\cap S)$
is a separating set with less than $k$ vertices since $|F\cap S_{1}|+|S_{1}\cap S|+|\bar{F}\cap S_{1}|=|S_{1}|=k$,
a contradiction. For another, if $\bar{F}\cap \bar{F}_{1}\neq \emptyset$, it yields that
$(\bar{F}\cap S_{1})\cup (S_{1}\cap S)\cup (\bar{F}_{1}\cap S)$ is a separating set with less than $k$ vertices, also a contradiction.
So $F\cap \bar{F}_{1}= \emptyset$ and $\bar{F}\cap \bar{F}_{1}= \emptyset$.
Then we again by the fact that $\bar{F}_{1}\subset S$ and $\delta(\bar{F}_{1})\geq \delta(G)-|S_{1}|\geq 1$  conclude a contradiction as above way.

\begin{case}
$|S_{1}\cap F|>|\bar{F}\cap S_{1}|$ or $|S_{1}\cap F|<|\bar{F}\cap S_{1}|$.
\end{case}
Suppose that $|S_{1}\cap F|>|\bar{F}\cap S_{1}|$ (For $|S_{1}\cap F|<|\bar{F}\cap S_{1}|$, it can be shown by the same proof way, we omit it).
Without loss of generality, let $|F_{1}\cap S|\geq|\bar{F}_{1}\cap S|$. If $\bar{F}_{1}\cap \bar{F}\neq \emptyset$,
$(\bar{F}\cap S_{1})\cup (S_{1}\cap S) \cup (\bar{F}_{1}\cap S)$ is a separating set with less than $k$ vertices, a contradiction. So $\bar{F}_{1}\cap \bar{F}= \emptyset$.
If $F_{1}\cap \bar{F}=\emptyset$, it follows that $\bar{F}\cap S_{1}=\bar{F}\neq \emptyset$.
Since $\delta(\bar{F})\geq \delta(G)-|S|\geq 1$, it follows that $\Bar{F}$ contains at least an edge, say $x'y'$. Owing to $\delta_{G}(X')\geq k+1$ and $\delta_{G}(Y')\geq k+1$, we have $d_{G}(x')+d_{G}(y')\geq 2(k+1)$. Meanwhile, since $x'$ and $y'$ have distinct neighbors in $S$, we also have $d_{G}(x')+d_{G}(y') \leq |S|+|\Bar{F}|\leq |S|+|S_{1}|\leq 2k$, a contradiction.
It results in $F_{1}\cap \bar{F}\neq\emptyset$.
But if $\bar{F}_{1}\cap F=\emptyset$, then $\bar{F}_{1}\cap S=\bar{F}_{1}\neq \emptyset$.

Since $\delta(\bar{F}_{1})\geq \delta(G)-|S_{1}|\geq 1$,  $\Bar{F_{1}}$ contains at least an edge, say $x''y''$. Owing to $\delta_{G}(X')\geq k+1$ and $\delta_{G}(Y')\geq k+1$, we have $d_{G}(x'')+d_{G}(y'')\geq 2(k+1)$. Meanwhile, since $x''$ and $y''$ have distinct neighbors in $S$, we also have $d_{G}(x'')+d_{G}(y'') \leq |S_{1}|+|\Bar{F_{1}}|\leq |S_{1}|+|S|\leq 2k$, a contradiction.
Then it contributes to $\bar{F}_{1}\cap F\neq\emptyset$.

As above analysis, we have that $\bar{F}_{1}\cap \bar{F}= \emptyset$, $\bar{F}_{1}\cap F\neq\emptyset$,
and $F_{1}\cap \bar{F}\neq\emptyset$. Now we claim that $|F\cap S_{1}|>|S\cap F_{1}|$. 
Otherwise, if $|F\cap S_{1}|=|S\cap F_{1}|$, then $\bar{F_{1}}\cap F$ is a fragment which contradicts the fact that $F$ is an end; if $|F\cap S_{1}|<|S\cap F_{1}|$, then
the set $(F\cap S_{1})\cup(S\cap S_{1})\cup (\bar{F}_{1}\cap S)$ is a separating set with size less than $k$ vertices,
$(F\cap S_{1})\cup(S\cap S_{1})\cup (\bar{F}_{1}\cap S)$ is a separating set with less than $k$ vertices,
which contradicts the fact that $F$ is an end.
According to $|F\cap S_{1}|+|S_{1}\cap S|+|\bar{F}\cap S_{1}|=|S_{1}|=k$ and $|F\cap S_{1}|>|S\cap F_{1}|$, it follows that $(\bar{F}\cap S_{1})\cup (S_{1}\cap S)\cup (F_{1}\cap S)$ is a separating set with
less than $k$, a contradiction.
This completes the proof.
\end{proof}

Before proceeding the proof, we reminder of something that if deleting two adjacent vertices of $G$, the minimum degree of $G$ decreases at most by one, because $G$ is a bipartite graph. Let $S$ be a minimum separating set of $G$ and $F$ an end to $S$. Then we consider
the cases that $\kappa(G)=k$ and $\kappa(G)>k$, respectively.
\setcounter{case}{0}
\begin{case}\label{case1}
$\kappa(G)=k$.
\end{case}
If $m=1$, then $\delta(G)\geq k+1$, 
and by Lemma \ref{lem1}, any vertex in $F$, say $v_{1}$, is not contained in any minimum separating set of $G$. Set $G_{1}:=G-\{v_{1}\}$. It implies that $G_{1}$ must be $k$-connected. So $G$ contains a path $P$ on one vertex.

Suppose $m\geq 3$.
If $F-\{v_{1}\}$ is an end of $G_{1}$ to $S$,
then, by Lemma \ref{lem1},
we continue to choose one neighbor of $v_{1}$, say $v_{2}\in (F-\{v_{1}\})$,
such that $G_{1}-\{v_{2}\}$ is $k$-connected.
If $F-\{v_{1}\}$ is not an end of $G_{1}$ to $S$,
then there is an end $F_{1}\subset F$ of $G_{1}$. 
Otherwise, $F$ contains another fragment, which contradicts the fact that $F$ is an end.
We say that there is a vertex $v_{2}\in F_{1}$ adjacent to $v_{1}$. Otherwise, $F_{1}$ will be an end of $G$, which contradicts the fact that $F$ is an end of $G$.
Thus, by Lemma \ref{lem1},
$G-\{v_{1},v_{2}\}=G_{1}-v_{2}$ is $k$-connected.
Let $G_{2}:=G_{1}-v_{2}$, then $\delta(G_{2})\geq k+\lceil \frac{m}{2}\rceil-1$. Now we can continuously delete vertices of $G_{2}$ as above applying Lemma \ref{lem1} $m-2$ times. Consequently, there is a path $P=v_{1}v_{2}\cdots v_{m}$ such that $\kappa(G-V(P))=k$.   

\begin{case}
$\kappa(G)>k$.
\end{case}
Let $\kappa(G)=k'$, and let $F$ be an end of $G$ to a minimum separating set $S'$ of size $k'$.
Let $P_{s}=v_{1}v_{2}\cdots v_{s}$ be a path
such that $\kappa(G-P_{s-1})=k+1$ and $\kappa(G-P_{s})=k$, where $s\leq m$. 
Let $G_{0}:=G-P_{s}$.
Then we have that $\delta(G_{0})\geq k+\lceil\frac{m}{2}\rceil-\lceil\frac{s}{2}\rceil$. Let $F_{0}$
be an end of $G_{0}$ to a minimum separating set $S_{0}$ of size $k$. Clearly,
there is a vertex $u_{1}\in N(v_{s})\cap F_{0}$. Otherwise,
we have that $\kappa(G-(P_{s-1}))=k$, a contradiction. If $s=m$, then we find out a path $P$ on $m$ vertices such that $\kappa(G-V(P))\geq k$. If $s<m$, then we apply Lemma \ref{lem1} $2(\lceil\frac{m}{2}\rceil-\lceil\frac{s}{2}\rceil)$ times as the proof of Case \ref{case1}. This completes the proof.
\vskip 2mm

\noindent{\bf Acknowledgment}

The author would like to thank the referee for careful reading and comments improving the presentation of the results.

\vskip 2mm
\noindent{\bf Declaration of interests}

The authors declare that they have no known competing financial interests or personal relationships
that could have appeared to influence the work reported in this paper.


\begin{thebibliography}{1}

\bibitem{BM}
J. A. Bondy and U. S. R. Murty, \emph{Graph Theory}, Graduate Texts
in Mathematics 244, Springer, Berlin, 2008.

\bibitem{CKL}
G. Chartrand, A. Kaigars and D.R. Lick, Critically $n$-connected
graphs, \emph{Proc. Amer. Math. Soc.} 32 (1972), 63--68.


\bibitem{AN}
A.A. Diwan and N.P. Tholiya, Non-separating trees in connected graphs,
\emph{Discrete Math.} 309 (2009), 5235--5237.

\bibitem{SK}
S. Fujita and K. Kawarabayashi, Connectivity keeping edges in graphs
with large minimum degree, \emph{J. Combin. Theory, Ser. B} 98
(2008), 805--811.

\bibitem{H}
T. Hasunuma, Connectivity keeping trees in $2$-connected graphs with
girth conditions, \emph{Combin. Algor.} (2020), 316--329.

\bibitem{TK}
T. Hasunuma and K. Ono, Connectivity keeping trees in $2$-connected
graphs, \emph{J. Graph Theory} 94 (2020), 20--29.

\bibitem{HL}
Y. Hong and Q. Liu, Mader's conjecture for graphs with small connectivity,  \emph{J. Graph Theory} 101(3) 2022, 379--388.

\bibitem{HL2}
Y. Hong, Q. Liu, C. Lu and Q. Ye, Connectivity keeping caterpillars and spiders in 2-connected graphs,
\emph{Discrete Math.} 344(3) (2021), 112236.

\bibitem{Lu}
C. Lv and P. Zhang, Connectivity keeping trees in $2$-connected graphs,
\emph{Discrete Math.} 343 (2) (2020), 1--4.

\bibitem{LTW}
L. Luo, Y. Tian and L. Wu, Connectivity keeping paths in $k$-connected bipartite graphs,
Discrete Math. 345 (4) (2022) 112788.

\bibitem{mader1}
W. Mader, Connectivity keeping paths in $k$-connected graphs,
\emph{J. Graph Theory} 65 (2010), 61--69.

\bibitem{mader2}
W. Mader, Connectivity keeping trees in $k$-connected graphs,
\emph{J. Graph Theory} 69 (2012), 324--329.

\bibitem{Tian1}
Y. Tian, H. Lai, L. Xu and J. Meng, Nonseparating trees in
$2$-connected graphs and oriented trees in strongly connected
digraphs, \emph{Discrete Math.} 342 (2019), 344--351.

\bibitem{Tian2}
Y. Tian, J. Meng and L. Xu, Connectivity keeping stars or double-stars
in $2$-connected graphs, \emph{Discrete Math.} 341 (2018),
1120--1124.

\bibitem{YT}
Q. Yang and Y. Tian, Connectivity keeping caterpillars and spiders in bipartite graphs with
connectivity at most three, \emph{Discrete Math.} 346 (2023), 113207.

\bibitem{YT2}
Q. Yang and Y. Tian, Proof of a conjecture on connectivity keeping odd paths in
k-connected bipartite graphs, arXiv:2209.08373v2.

\bibitem{Z}
P. Zhang, Research on the Existence of Subgraphs Related to Connectivity (Ph.D. dissertation), East China Normal University, 2021.

\end{thebibliography}
\end{document}